\documentclass[a4paper,12pt]{amsart}
\usepackage{amsfonts,amsthm,amssymb}
\usepackage{epsfig}
\usepackage{pstricks,slashbox,multirow}
\newtheorem{theorem}{Theorem}[section]

\newtheorem{cor}[theorem]{Corollary}
\newtheorem{dnt}[theorem]{Definition}

\newtheorem{exm}[theorem]{Example}

\newtheorem{rem}[theorem]{Remark}

\title[Constant Sum Partition of $I_n$ With Prescribed Subsets Orders]{Constant Sum Partition of $\{1,2,...,n\}$ Into Subsets With Prescribed Orders}

\author{\sc V. Vilfred $Kamalappan^{1}$ ~and ~Sajidha ~$P^{2}$}
\address{Department of Mathematics, Central University of Kerala, Periye \\
Kasaragod, Kerala,  India - 671 320}
\email{$^1$vilfredkamal@gmail.com, $^2$sajinada555@gmail.com}

\date{}

\begin{document}
	
	\maketitle
	
	\begin{abstract}
Studies on partition of $I_n$ = $\{1, 2, . . . , n\}$ into subsets $S_1, S_2, . . . , S_x$ so far considered with prescribed sum of the elements in each subset. In this paper, we study constant sum partitions $\{S_1,S_2,...,S_x\}$ of $I_n$ with prescribed $|S_i|$, $1 \leq i \leq x$. Theorem \ref{thm 2.3} is the main result which gives a necessary and sufficient condition for a partition set $\{S_1,S_2,\ldots, S_x\}$ of $I_n$ with prescribed $|S_i|$ to be a constant sum partition of $I_n$, $1 \leq i \leq x$ and $n > x \geq 2$. We state its  applications in graph theory and also define {\em constant sum partition permutation} or {\em magic partition permutation} of $I_n$. A partition $\{S_1,S_2,\cdots,S_x\}$ of $I_n$ is a {\em constant sum partition of $I_n$} if $\sum_{j\in S_i}{j}$ is a constant for every $i$, $1 \leq i \leq x$. 
\end{abstract}

\section{Introduction}

Partition of numbers plays an important role in Combinatorics, Representation theory, Lie theory, mathematical physics, and theory of special functions \cite{gk}. Euler, Ramajuan, Rademacher, Paul Erodes and many others studied partitioning of the set $I_n$ = $\{1,2,...,n\}$ under different conditions and revealed its beauty and applications \cite{agk, cg, llwz, sa, ss}. A {\em partition of a set} is a grouping of its elements into non-empty subsets, in such a way that every element is included in exactly one subset.	{\em Partitioning of a natural number} $n$ is a way of writing $n$ as a sum of positive integers where the order of the addends is not significant, possibly subject to one or more additional constraints. Partitioning is used to solve maths problems involving large numbers easier by separating them into smaller units and constant sum partitions are used to construct magic rectangles which are a natural generalization of magic squares (\cite{rb}, pages 215-221). $k$-distance magic ($k$-DM) labeling \cite{v22} is a generalisation of DM labeling \cite{mrs,sf}  or sigma labeling \cite{v96}. 

In  \cite{llwz}, a nonincreasing sequence of positive integers $\left < m_1, m_2, . . . , m_x \right >$ is called  {\em $n$-realizable} if the set $I_n$ can be partitioned into $S_1, S_2, . . . , S_x$ such that $\sum_{j\in S_i}{j}$ = $m_i$ for each $i$, $1 \leq i \leq x$. And it is proved that nonincreasing sequence of positive integers $\left < m_1, m_2, . . . , m_x \right >$ is $n$-realizable under the conditions that $\sum^{x}_{i=1}{m_i}$ = $n(n+1)/2$	and $m_{x-1} \geq n$. In  \cite{ss}, partition $\{S_1,S_2,...,S_x\}$ of $I_n$ is called an {\em $(x,t)$-partition of $I_n$} when $m_i$ = $t$ for each $i$, $1 \leq i \leq x$.  And existence of $(x,t)$-partition of $I_n$ is established whenever $n(n+1)$ = $2xt$ and also an $r$-regular path-perfect graph $G$ with $p$ vertices and $n(n+1)/2$ edges using the existence of an $(r/2, p)$-partition of $I_n$ when $r$ is even. 

\begin{dnt} [\rm {\bf 1.1}] \cite{llwz} For $ n > x > 1$, a partition $ S_1,S_2,\ldots, S_x $ of $I_n$ is called a \emph{constant sum partition} if $\sum_{a \in S_i}a = M,~ a~ constant,$ $\forall i$, $1 \leq i \leq x$.  
\end{dnt}
Here, constant $M$ is referred as {\em Partite Constant} or {\em Magic Constant}.

This paper contains four sections. Section 1 is an introductory one. Section 2 contains our main result: Let $n > x \geq 2$, $n,x\in\mathbb{N}$, $\mathfrak{S}$ = $\{S_1,S_2,...,S_x\}$ be a partition of $I_n$ with prescribed $|S_i| = n_i$, $1 \leq n_1 \leq n_2 \leq . . . \leq n_x < n$, $N_i$ = $n_1+...+n_i$, $M$ = $\frac{T_{n}}{x}$ and $ 1 \leq i \leq x$. Then, $\mathfrak{S}$ is a constant sum partition of $I_n$ if and only if $M = \frac{T_{n}}{x}\in\mathbb{N}$ and $iM ~ \leq  \sum_{j =1}^{N_i}{(n-j+1})$ for all $i$, $ 1 \leq i \leq x$. In Section 3, we define {\em constant sum partition permutation} or {\em magic partition permutation} of $I_n$ and illustrate the method involved in the proof of the main result with examples. Applications of constant sum partitions of $I_n$ in graph theory is presented in Section 4. 

Effort to study $k$-DM graphs is the motivation for this work. 

\section{ \bf On $x$-Partite Constant Sum Partition of $I_{n}$}

Here, we study constant sum partition  $\{S_1,S_2,\cdots,S_x\}$ of $I_n$ with prescribed $|S_i|$ for all $i$, $1 \leq i \leq x < n$ and $ n \geq 3$. Theorem \ref{thm 2.3} is the main result which gives a necessary and sufficient condition for a partition set $\{S_1,S_2,\ldots, S_x\}$ of $I_n$ with prescribed $|S_i|$ to be a constant sum partition of $I_n$, $1 \leq i \leq x$ and $n > x \geq 2$. 

\begin{dnt} [\rm {\bf 2.1}] \cite{db} \quad  A \emph{triangular number} is a natural number which is the sum of finite consecutive positive integers, beginning with 1. i.e., a number is \emph{triangular} if it is of the form $ \frac{n(n+1)}{2} $ for some $n\in\mathbb{N}.$ 
	
	Triangular number corresponding to the sum of the first $n$ consecutive natural numbers is denoted by $T_n$. i.e., $T_n$ = $\frac{n(n+1)}{2}$, $n,T_n\in\mathbb{N}.$
\end{dnt}

Let $ n > x > 1$ and $\{ S_1,S_2,\ldots, S_x\} $ be an $x$-partite constant sum partition of $I_n$ with magic constant $M$. Then, 
$$\sum_{i =1}^{n}{i} =  \sum_{a\in S_1}a +  \sum_{a\in S_2}a + \cdots +  \sum_{a\in S_x}a.~$$ 

$$\Rightarrow~	\sum_{i =1}^{n}{i} =  x \sum_{a\in S_i}a,~~~ 1 \leq i \leq x.  \hspace{2.5cm} $$ 

$$\Rightarrow~~~T_n = M x.~ \Rightarrow~ T_{n} \equiv 0~ (mod~ x), ~~M,T_n,x\in\mathbb{N}. \hspace{1.2cm} (1)$$
Thereby, we get the following result.

	\begin{theorem}\label{thm b1} \quad {\rm Let $ n > x > 1$ and $M$ be the magic constant of $x$-partite constant sum partition of $I_n$, $ n,x,M \in \mathbb{N}$. Then, $M$ = $\frac{T_n}{x}\in\mathbb{N}$ and $ T_{n} \equiv 0~ (mod~ x)$. \hfill $\Box$}
\end{theorem}

\begin{theorem}\label{thm 2.3} \quad {\rm 
		Let $n > x \geq 2$, $n,x\in\mathbb{N}$, $\mathfrak{S}$ = $\{S_1,S_2,...,S_x\}$ be a partition of $I_n$ with prescribed $|S_i| = n_i$, $N_i$ = $n_1+...+n_i$, $1 \leq n_1 \leq n_2 \leq . . . \leq n_x < n$, $M$ = $\frac{T_{n}}{x}$ and $ 1 \leq i \leq x$. Then, $\mathfrak{S}$ is an $x$-partite constant sum partition of $I_n$ if and only if $M = \frac{T_{n}}{x}\in\mathbb{N}$ and $iM ~ \leq  \sum_{j =1}^{N_i}{(n-j+1})$ = $T_n -  T_{n-N_i}$ for all $i$, $ 1 \leq i \leq x$. }
\end{theorem}
\begin{proof}\quad 	Let $ \{S_1,S_2,\ldots,S_x\} $ be a constant sum partition of $ I_{n} $ with magic constant $M$ = $M_i$ = $\displaystyle \sum_{a\in S_i}a$, $|S_i| = n_i$ and $N_i$ = $n_1+...+n_i$, $1 \leq i \leq x$. Clearly, for all $i$ $\ni$ $ 1 \leq i \leq x $, $M = \frac{T_{n}}{x}\in\mathbb{N}$ and 
	
	$\displaystyle 	\sum_{j =1}^{N_i}{(n-j+1})$ = $\displaystyle 
	\sum_{j =1}^{n-(n-N_i)}{(n-j+1})$ 
	
\hspace{2.5cm}	= $\displaystyle 
	\sum_{j =1}^{n}{(n-j+1})$ -  $\displaystyle 
	\sum_{j =1}^{n-N_i}{(n-j+1})$ = $T_n -  T_{n-N_i}$.  \\
	Also, for $ ~1 \leq i \leq x,$ sets $ S_i $'s are mutually disjoint and 
	$$ \sum_{a\in {S_1 \cup S_2 \cup \cdots  \cup  S_i}}{a} ~~\leq \sum_{j =1}^{n_1+n_2+\cdots + n_i}{(n-j+1)}. \hspace{4cm} $$
	$$	\Rightarrow	iM \leq \sum_{j =1}^{N_i}{(n-j+1)} = T_n -  T_{n-N_i},~ 1 \leq i \leq x. $$  
	
	Conversely, we prove that there exists an $x$-partite constant sum partition to $I_{n}$ under the given conditions. 
	
Given, $M = \frac{T_{n}}{x}\in\mathbb{N}$ implies, $Mx = T_{n} = \sum_{i=1}^n{i}\in\mathbb{N}$. Also, if for any $j$ $\ni$ $1 \leq j \leq x$, 	
	
	\hspace{2cm} $	jM ~ > \sum_{i =1}^{N_j}{n+1-i} = T_n -  T_{n-N_j},$
	\\
	then $x$-partite constant sum partition of $I_n$ doesn't exist since sum of elements of partite sets $S_1, S_2, . . . , S_j$ of $I_n$ must be $\leq$ sum of $N_{j}$ biggest numbers of $I_n$ for any $j$ when $\{S_{1}, S_{2}, . . . , S_{j}\}$ is a constant sum partition of $I_{n}$, $N_{j} = n_{1}+n_{2}+\cdots+n_{j}$, $n_{j}$ = $|S_{j}|$, $1 \leq n_1 \leq n_2 \leq$ $. . .$ $\leq n_x < n$ and $n_1,n_2,...,n_j\in\mathbb{N}$.
	
	Thus, for the partition set $\{S_1,S_2,\ldots,S_x\}$ of $I_n$ with prescribed $|S_i|$ = $n_i$ and $1 \leq n_1 \leq n_2 \leq . . . \leq n_x < n$ to be a constant sum partition set of $I_n$, $M = \frac{T_{n}}{x}\in\mathbb{N}$ and 
	$$	iM ~ \leq \sum_{j =1}^{N_i}{n+1-j} = T_n -  T_{n-N_i},~ \forall~ i, ~1 \leq i \leq x < n.$$
	Now, we are left with obtaining of an $x$-partite constant sum partition to $I_n$.
	
	Let $\{S_1,S_2,\ldots,S_x\}$ be a partition of $I_{n}$ with $|S_i|$ = $n_i$, $n > x > 1$, $M$ = $\frac{T_n}{x}$, $ 1 \leq n_1 < n_2 \leq \cdots \leq n_x$, $N_i = n_1 + n_2 +\cdots + n_i $, $ n = N_x,$ $n,x,M,N_i \in \mathbb{N}$ and 
	$$ 	iM ~ \leq \sum_{j =1}^{N_i}{n+1-j} = T_n -  T_{n-N_i} ~\forall~ i, ~ 1 \leq i \leq x.$$
	$$\Rightarrow~	M ~ \leq \frac{1}{i} \sum_{j =1}^{N_i}{(n-j+1}),  ~\forall~ i, ~ 1 \leq i \leq x. \hspace{2cm} (2)$$
	
	If $\{S_1,S_2,\ldots,S_x\}$ is a constant sum partition of $I_{n}$ with $M$ = $M_i$ = $\sum_{j\in S_i}j$ for $i$ = 1 to $x$, then the proof is over. Otherwise, consider the following algorithm to obtain an $x$-partite constant sum partition $\{S_1,S_2,\ldots,S_x\}$ of $I_{n}$.
	
	\vspace{.15cm}	 
	\noindent \textbf{Algorithm:}\\
	\noindent \textbf{\textit{Step 1.}}
	Arrange elements of $I_n$ in descending order as $I_n = \{n,n-1,$ $n-2,$ $\ldots,1\}$.	
	
	\vspace{.15cm}	 
	\noindent \textbf{\textit{Step 2.}}
	Let $A_1$ = $\left < a_{1,1},a_{1,2},\ldots,a_{1,n_1}\right >$ and $W_1 =  \sum_{j=1}^{n_1}{a_{1,j}}$ where $a_{1,j}$ is the $j^{th}$ element $n-j+1$ of $I_n$, $1 \leq j \leq n_1$. 
	
	$\Rightarrow$ $A_1$ = $\left < n,n-1,\ldots,n-(n_1-1)\right >$.
	
	Using (2), $W_1 \geq M$. We consider the following two cases.
	
	\vspace{.15cm}	 
	In the first case, $  M = W_1 = \sum_{a\in A_1}{a}$. Take $S_1$ = $\{a_{1,1},a_{1,2},...,a_{1,n_1}\}$.
	
	In the second case, $W_1 > M $ which implies, $W_1 - M \in \mathbb{N}$. 
	In this case, let $B_1$ = $< 0,0,\ldots, 0, p_{1,1},  p_{1,2}, . . . ,$ $p_{1,k_1} >$ such that $|B_1|$ = $n_1$, $p_{1,1} + p_{1,2} + . . . + p_{1,k_1}$ = $W_1-M$, $p_{1,1} \leq p_{1,2} \leq . . .$ $\leq p_{1,k_1} < a_{1,n_1}$ and $p_{1,j}\in\mathbb{N}$, $1 \leq j \leq k_1 \leq n_1$. 
	
	$\therefore$ $A_1 - B_1$ = $< a_{1,1},a_{1,2},\ldots,a_{1,n_1-k_1},  a_{1,n_1-(k_1-1)}- p_{1,1}, a_{1,n_1-(k_1-2)}- p_{1,2}, . . . ,$ $a_{1,n_1}-p_{1,k_1} >$. If any pair of elements of $A_1 - B_1$ are not distinct, then choose a different partition of $W_1-M$ such that all the elements of $A_1 - B_1$ are distinct and $A_1 - B_1 \subseteq I_n$.  Existence of such set $A_1 - B_1$ is possible because of the conditions on $p_{1,j}$ and of relation $(2)$, $1 \leq j \leq k_1 \leq n_1$.  Clearly,  
	$$ \sum_{a\in A_1 - B_1} a = \sum_{a\in A_1} a - \sum_{a\in B_1} a = \sum_{j=1}^{n_1} a_{1,j} -  \sum_{j=1}^{k_1} p_{1,j} = W_1 - (W_1 -M) = M.$$
	Now, take $S_1$ = $\{a_{1,1},a_{1,2},\ldots,a_{1,n_1-k_1},  a_{1,n_1-(k_1-1)}- p_{1,1}, a_{1,n_1-(k_1-2)}- p_{1,2}, . . . ,$ $a_{1,n_1}-p_{1,k_1}\}$. After getting $S_1$, we find $S_2$, . . . , $S_x$ one by one by the following general method. 

\vspace{.15cm}	 
\noindent \textbf{Step 3.}
We find $S_i$ after obtaining $S_1$, . . . , $S_{i-1}$ as follows, $2 \leq i \leq x$. 
	
	Let $A_i$ = $\left < a_{i,1},a_{i,2},\ldots,a_{i,n_i}\right >$ = $<$ elements of $I_n \setminus (S_1 \cup . . . \cup S_{i-1})$ as components arranged in descending order $>$, $2 \leq i \leq x$. Let $W_i = \displaystyle \sum_{j =1}^{n_i}{a_{i,j}}$, $1 \leq i \leq x$. For $j$ = 1 to $i-1$, $W_j$ = $M$, $2 \leq i \leq x <n$. Using (2), $W_i \geq M$, $2 \leq i \leq x$. 
	
	\vspace{.1cm}			 
	When $  W_i = M$,  take $S_i$ = $\{a_{i,1},a_{i,2},...,a_{i,n_i}\}$, $2 \leq i \leq x$.
	
	\vspace{.1cm}			 
	When $W_i > M $, $W_i - M \in \mathbb{N}$, $2 \leq i \leq x$. 
	In this case, let $B_i$ = $<0,0,\ldots, 0,$ $p_{i,1}, p_{i,2}, . . . , p_{i,k_i}>$ such that $|B_i|$ = $n_i$, $p_{i,1}$ + $p_{i,2}$ + . . . + $p_{i,k_i}$ = $W_i-M$, $p_{i,1} \leq p_{i,2} \leq . . . \leq p_{i,k_i} < a_{i,n_i}$ and $p_{i,j}\in\mathbb{N}$, $1 \leq j \leq k_i \leq n_i$. 
	
	Let $A_i - B_i = ~ <a_{i,1},a_{i,2},\ldots,a_{i,n_i-k_i},  a_{i,n_i-(k_i-1)}- p_{i,1}, a_{i,n_i-(k_i-2)}- p_{i,2}, . . . ,$ $a_{i,n_i}-p_{i,k_i}>$. If any pair of elements/components of $A_i - B_i$ are not distinct, then choose a different partition of $W_i-M$ such that all the elements/components of $A_i - B_i$ are distinct and belong to $I_n \setminus(S_1 \cup S_2 \cup . . . \cup S_{i-1})$, $2 \leq i \leq x$. Clearly,  
	
	$ \sum_{a\in A_i - B_i} a = \sum_{a\in A_i} a - \sum_{a\in B_i} a = \sum_{j=1}^{n_i} a_{i,j} -  \sum_{j=1}^{k_i} p_{i,j}$ 
	
	\hfill $ = W_i - (W_i -M) = M.$
	\\
	In the above, for example, for a given $i$, one can choose $p_{i,j}$ and $m_i$ such that for $j$ = 1 to $m_i$, $p_{i,j}$ = $j$,  $\sum_{j=1}^{m_i-1}{p_{i,j}} < W_i-M \leq \sum_{j=1}^{m_i}{p_{i,j}}$, $m_i < n_i$ and $2 \leq i \leq x$. And thereafter one has to adjust the elements $p_{i,j}$ either by removal, replacement or by both but still the resultant elements of $A_i-B_i$ should be distinct, belong to $I_n \setminus(S_1 \cup S_2 \cup . . . \cup S_{i-1})$ and their sum equal to $M$. At this stage, take $S_i$ = $\{a_{i,1}, a_{i,2}, . . . , a_{i,n_i-k_i},  a_{i,n_i-(k_i-1)} - p_{i,1},$ $a_{i,n_i-(k_i-2)} - p_{i,2}, . . . , a_{i,n_i-(k_i-k_i)}-p_{i,k_i}\}$ = $\{a_{i,1}, a_{i,2}, . . . , a_{i,n_i-k_i},  a_{i,n_i-(k_i-1)} - p_{i,1},$ $a_{i,n_i-(k_i-2)} - p_{i,2}, . . . , a_{i,n_i}-p_{i,k_i}\}$, $2 \leq i \leq x$. 
		
	\vspace{.1cm}			 
	\noindent \textbf{\textit{Step 4.}} 
	Continue Step 3 untill $S_x$ is obtained. 
	
	\vspace{.1cm}			 
	The above algorithm is well defined and the iterative  process will endup at finite iterations since $I_n$ is a finite set and $S_i$s are subsets of $I_n$, $1 \leq i \leq x < n$. Note that after obtaining $S_1, S_2, . . . , S_{x-1}$, the set $S_x$ = $I_n\setminus (S_1 \cup S_2 \cup . . . \cup S_{x-1})$.
	
	Thus, we could find an $x$-partite constant sum partition $\{S_1,S_2,\ldots,$ $S_x\}$ of $I_{n}$ with prescribed $|S_i|$ = $n_i$, $1 \leq n_1 \leq n_2 \leq . . . \leq n_x < n$, $n_1 + n_2 + ... + n_x$ = $n$, $n > x > 1$ and $M$ = $\displaystyle \sum_{a\in S_i}a$, $1 \leq i \leq x.$ Hence the result.
\end{proof}

Method of proof given in the above theorem is illustrated by Example \ref{e4} in finding constant sum partition $\{S_1,...,S_x\}$ of $I_{n}$ with $n$ = 105, $x$ = 5, $n_1 = 12 = n_2$, $n_3$ = 15, $n_4$ = 20 and $n_5$ = 46 where $n_i$ = $|S_i|$ and $1 \leq i \leq x < n$.  

\begin{cor}\label{thm 2.5} \quad {\rm 
		Let $n > x \geq p \geq 2$, $\mathfrak{S}$ = $\{S_1,S_2,\ldots,S_p,\ldots,S_x\}$ be a partition of $I_n$ with prescribed $|S_i| = n_i$, $ 1 \leq n_1 = n_2 = \cdots = n_p \leq n_{p+1} \leq ... \leq n_x$, $M$ = $\frac{T_{n}}{x}$, $ 1 \leq i \leq x$ and $n,p,x\in\mathbb{N}$. Then, $\mathfrak{S}$ is an $x$-partite constant sum partition of $I_n$ if and only if $M = \frac{T_n}{x}\in\mathbb{N}$, $M \leq\left \lfloor \frac{Q_q}{q} \right \rfloor \leq$ $\sum_{j =1}^{n_1}{(n-j+1)}$ and $iM \leq \sum_{j =1}^{n_1+n_2+\cdots + n_i}{(n-j+1)}$ for every $q$ and $i$ where $Q_q = \displaystyle  \sum_{j =1}^{qn_1}{(n-j+1)}$, $1 \leq q \leq p \leq x < n$ and $1 \leq i \leq x$.	}	
\end{cor}
\begin{proof} \quad For every $i$ and $q$ such that $n > x \geq p \geq q \geq 1$, $p \geq 2$ and $1 \leq i \leq x$,
	$$M_i =	\sum_{a\in S_i}a  \leq \sum_{j =1}^{n_1}{(n-j+1)}~ and~ \hspace{2cm}$$ 
	$$	\sum_{a\in {S_1 \cup S_2 \cup \cdots  \cup  S_q}}{a} \leq \sum_{j =1}^{n_1+n_2+\cdots + n_q}{(n-j+1)} = \sum_{j =1}^{qn_1}{(n-j+1)} = Q_q \hspace{.75cm} (3)$$
	since in $I_n$, sum of any $qn_1$ distinct elements is less than or equal to sum of $qn_1$ biggest elements.
	
	When $\mathfrak{S}$ is a constant sum partition of $I_n$, $M_i = M = \frac{T_n}{x}\in\mathbb{N}$, $\forall$ $i$, $1 \leq i \leq x.$
	$$ And ~ (3)~ \Rightarrow qM_i = qM = q ~\times ~ \frac{T_n}{x} \leq  Q_q,  ~ 1 \leq i,q \leq p. $$ 		
	$$\Rightarrow M_i = M = \frac{T_n}{x} \leq  \left \lfloor \frac{\sum_{j =1}^{qn_1}{(n-j+1)}}{q} \right \rfloor = \left \lfloor \frac{Q_q}{q} \right \rfloor ~ and~ M = \frac{T_n}{x}\in\mathbb{N}. \hspace{1cm} $$ 
	$$ 	\hspace{.5cm} \Rightarrow M_i = M \leq  \sum_{j =1}^{n_1}{(n-j+1)} ~ ~ and~ M = \frac{T_n}{x}\in\mathbb{N},~ \forall q, \hspace{2cm} (3)$$
	$1 \leq q \leq p \leq x $ since $ \frac{\sum_{j =1}^{qn_1}{(n-j+1)}}{q} \leq \sum_{j =1}^{n_1}{(n-j+1)}$ for $p \geq q > 1$ and $M\in\mathbb{N}$. Then the result follows from Theorem \ref{thm 2.3}.
\end{proof}

\begin{theorem}\label{thm 2.6}
	\quad {\rm  Let $n > x \geq p \geq 2,$ $\mathfrak{S}$ = $ \{S_1,S_2,\ldots,S_p, . . . ,S_x\} $ be a partition of $ I_{n}$,  $n_i = |S_i|$, $M_i =	\sum_{a\in S_i}a$, $1 \leq n_1 = n_2 = \cdots = n_p \leq n_{p+1} \leq \cdots \leq n_x $, $Q_p = \displaystyle  \sum_{j =1}^{pn_1}{(n-j+1)} $, $M = \frac{T_{n}}{x}$ and $n,x,p,n_i\in\mathbb{N}$. If $M > \left \lfloor \frac{Q_p}{p} \right \rfloor $, then $\mathfrak{S}$ can not be a constant sum partition  of $I_n$.}
\end{theorem}
\begin{proof}\quad When $ \left \lfloor M \right \rfloor > M$, clearly $\mathfrak{S}$ can not be an $x$-partite  constant sum partition of $I_n$. Now, assume that $ \left \lfloor M \right \rfloor = M$. 
	
	Suppose, $I_n$ be an $x$-partite constant sum partition with magic constant $M$. Then, for $ 1 < p \leq x $, using Theorem \ref{thm 2.3},  
	
	\vspace{.2cm}
	$pM  \leq \sum_{j =1}^{n_1+n_2+\cdots+n_p}{(n-j+1)} = \sum_{j =1}^{pn_1}{(n-j+1)} = Q_p$  which implies,	
\\
	$M \leq \left \lfloor \frac{Q_p}{p} \right \rfloor$ which is a contradiction to the given condition. Hence the result.
\end{proof}

\section{Illustration by examples on the proof of Theorem \ref{thm 2.3}}

Theorem \ref{thm 2.3} is the important result in this paper. In this section, we illustrate the method of the proof of this theorem with examples. Here, we start with the definition of {\em constant sum partition permutation} or {\em magic partition permutation} of $I_n$.  

Let $A$ be a nonempty set. A {\em permutation} of $A$ is a bijection from $A$ to itself.

\begin{dnt} \quad Let $X$ be non-empty set, $|X| \geq 2$, $\mathfrak{X}$ = $\{X_1, X_2,$ $\ldots, X_x\}$ be a partition of $X$ and $\Pi$ be a permutation on the set $X$. Then, permutation $\Pi$ is called a {\em permutation on partition $\mathfrak{X}$} or {\em partition permutation} on  $\mathfrak{X}$ if $\Pi$ preserves the partition. That is, under permutation $\Pi$ on $X$, $\{\Pi(X_1), \Pi(X_2), . . . , \Pi(X_x)\}$ is a partition of $X$ whenever $\{X_1,X_2,\ldots,X_x\}$ is a partition of $X$.  	
\end{dnt}

\begin{dnt} \quad Let $n > x \geq 2$, $\mathfrak{S}$ = $\{S_1,S_2,\ldots,S_x\}$ be a partition of $I_n$ with prescribed partite sums (orders) and $\Pi$ be a permutation on $I_n$. Then, permutation $\Pi$ is called an {\em $x$-partite  permutation on $\mathfrak{S}$ with prescribed partite sums (orders)} if $\Pi$ preserves the partition and sum (number) of the elements in $\Pi(S_i)$ and in $S_i$ are the same for every $i$, $1 \leq i \leq x$.   	
\end{dnt}

\begin{dnt} \quad Let $X \subset \mathbb{N}$, $|X| \geq 2$, $\mathfrak{X}$ = $\{X_1, X_2, . . . , X_x\}$ be a constant sum partition of $X$ (with prescribed partite orders) and $\Pi$ be a permutation on $X$. Then, $\Pi$ is called a {\em constant sum partition $\mathfrak{X}$ permutation on $X$} or {\em magic partition $\mathfrak{X}$ permutation on $X$} if $\Pi$ preserves constant sum partition $\mathfrak{X}$ (with prescribed partite orders). That is, under permutation $\Pi$ of $X$, $\{\Pi(X_1), \Pi(X_2), . . . , \Pi(X_x)\}$ is a constant sum partition of $X$ whenever $\{X_1,X_2,\ldots,X_x\}$ is a constant sum partition of $X$ (with $|\Pi(X_i)|$ = $|X_i|$ for every $i$, $1 \leq i \leq x$).   	
\end{dnt}

\begin{dnt} \quad Let $n > x \geq 2$, $\mathfrak{S}$ = $\{S_1,S_2,...,S_x\}$ be an $x$-partite constant sum partition of $I_n$ (with prescribed partite orders) and $\Pi$ be a permutation on $I_n$. Then, $\Pi$ is called a {\em constant sum partition  $\mathfrak{S}$ permutation} or {\em magic partition $\mathfrak{S}$ permutation} of $I_n$ (with prescribed partite orders) if $\Pi$ preserves constant sum partition $\mathfrak{S}$ (with prescribed partite orders).	That is, under permutation $\Pi$ on $I_n$, $\{\Pi(S_1), \Pi(S_2), . . . ,$ $\Pi(S_x)\}$ is a constant sum partition of $I_n$ whenever $\{S_1,S_2,\ldots,S_x\}$ is a constant sum partition of $I_n$ (with $|\Pi(S_i)|$ = $|S_i|$ for every $i$, $1 \leq i \leq x$).   	
\end{dnt}

\begin{rem} \label{rem 5} \quad Thus, to obtain different $x$-partite constant sum partitions of $I_n$ w.r.t. a given $x$-partite constant sum partition $\mathfrak{S}$ of $I_n$ (with prescribed partite orders), it is enough to find all possible $x$-partite magic $\mathfrak{S}$ permutations of $I_n$. 
\end{rem}

Examples \ref{e1} to \ref{e4} show uses of Theorems \ref{thm 2.3} to \ref{thm 2.6} and Example \ref{e4} illustrates the method of proof of Theorem \ref{thm 2.3}.

\begin{exm} \quad \label{e1}
	4-partite constant sum partition of $I_{10}$ doesn't exist.
	
	{\rm Here, $n$ = 10, $x$ = 4 and $ T_n = T_{10} =   \frac{10\times 11}{2}$. This implies that $ M = \frac{T_{10}}{4} \notin\mathbb{N}$ and thereby 4-partite constant sum partition  of $ I_{10}$ doesn't exist by Theorem \ref{thm 2.3}. \hfill $\Box$}
\end{exm}
\begin{exm}\quad \label{e2}
	Find constant sum partitions of $I_{10}$, if exist.  
	
	{\rm  Here, $n$ = 10 and $T_n$ = $T_{10} = 55$. 5 and 11 are the only proper divisors of $T_{10}$.  This implies, 5 and 11 are the possible values of $x$ for which $I_{10}$ has $x$-partite constant sum partition so that $M$ = $T_n/x\in\mathbb{N}$. In the case of $x$ = 11, $x = 11 > 10 = n$ and hence 11-partite constant sum partition of $I_{10}$ doesn't exist. Now, consider the  case $x$ = 5 for possible 5-partite constant sum partition of $I_{10}$. 
		
		If possible, let $\{S_1, S_2, S_3, S_4, S_5\}$ be a 5-partite constant sum partition of $I_{10}$, $|S_i|= n_i$, $ 1 \leq n_1 \leq n_2 \leq n_3 \leq n_4 \leq n_5 $, $M$ = $T_{10}/5$ = 11 and  $ n_1 + n_2 + n_3 + n_4 + n_5 = 10$, $1 \leq i \leq 5$. Then different possible set of values of $\{ n_1, n_2, n_3, n_4, n_5\}$ are
			(i) $\{1, 2, 2, 2, 3\}$ and (ii) $\{2, 2, 2, 2, 2\}$. 
		
		Case (i) of $\{ n_1, n_2, n_3, n_4, n_5\}$ = $\{1, 2, 2, 2, 3\}$ is not possible since $M$ = $11 > 10$ = $n$, the biggest element in $I_{10}$ and the possible element of $S_1$ nd the result follows by Theorem \ref{thm 2.6}.
		
		In case (ii), for $S_1$ = $\{10,1\}$, $S_2$ = $\{9,2\}$, $S_3$ = $\{8,3\}$, $S_4$ = $\{7,4\}$ and $S_5$ = $\{6,5\}$, $\{S_1, S_2, S_3, S_4, S_5\}$ is a 5-partite constant sum partition of $I_{10}$. 
	Thus, $\{ \{10, 1\},$ $\{9, 2\},$ $\{8, 3\},$ $\{7, 4\},$  $\{6, 5\}\}$ is a (the only 5-partite) constant sum partition of $I_{10}$. \hfill $\Box$}
\end{exm}

\begin{exm} \quad \label{e3}
	Find a constant sum partition $ \{S_1, S_2, S_3\} $ of $ I_{20} $, if exists, where $|S_i| = n_i$, $ n_1$ = 4 = $n_2$ and $n_3$ = 12.
	
	\vspace{.1cm}
	\noindent
	{\rm 	Here, $ x$ = 3,  $n_1 = n_2$ = 4 and  $n_3 = 12$ which implies, $p = 2$, $ M = \frac{T_{20}}{3} = 70$ and 
		$$\frac{1}{p}\sum_{i =1}^{pn_1}{(n-i+1)} = \frac{1}{2} \sum_{i =1}^{2n_1}{(n-i+1)} = 66.$$
		$$\Rightarrow M = 70 > 66 = \frac{1}{p}\sum_{i =1}^{pn_1}{(n-i+1)}.$$  
		$\therefore$ Using Theorem \ref{thm 2.3}, $3$-partite constant sum partition of $I_{20}$ does not exist.  \hfill $\Box$
	}
\end{exm}

\begin{exm} \quad \label{e4} Find a constant sum partition $\{S_1, S_2, S_3, S_4, S_5\} $ of $ I_{105}$, if it exists, with $n_1$ = 12 = $n_2$, $n_3$ = 15, $n_4$ = 20, $n_5$ = 46 and $|S_i | = n_i$, $1 \leq i \leq 5$.
	
	\vspace{.2cm}
	\noindent
	{\rm Here, $n$ = 105, $x$ = 5, $M$ = $T_n/x$ = $\frac{105\times 106}{2\times 5}$ = $1113\in\mathbb{N}$ and
	\\
	 $W_1 = \sum_{i =1}^{n_1}{(n-i+1)} = \sum_{i =1}^{12}{(105-i+1)}$
	 
	 \hfill $= 105+...+94 = 1194 > 1113 = M;$ 
		\\
		$\sum_{i =1}^{n_1+n_2}{(n-i+1)} = \sum_{i =1}^{24}{(105-i+1)}$
		
		\hfill $ = 105+...+82 = 2244 > 2226 = 2M;$
		\\
		$\sum_{i =1}^{n_1+n_2+n_3}{(n-i+1)} = \sum_{i =1}^{39}{(105-i+1)}$
		
		\hfill $ = 105+...+67 = 3354 > 3M;$ 
		\\
		$\sum_{i =1}^{n_1+n_2+n_3+n_4}{(n-i+1)} = \sum_{i =1}^{59}{(105-i+1)}$
		
		\hfill $ = 105+...+47 = 4484 > 4M.$
		
		Thus, $iM \leq \sum_{i=1}^{n_1+...+n_i}{(n-i+1)}$ for all  $i$, $1 \leq i \leq x = 5$ and thereby we use the method given in Theorem \ref{thm 2.3} to obtain 5-partite constant sum partition of $I_{105}$ with given $n_j$, $j$ = 1 to 5. 
		
		\vspace{.2cm}		
		\noindent	
		{\it Finding the 5-partite constant sum partition of $I_{105}$:}
		
		\vspace{.1cm}
		\noindent		
		\item [\rm (i)]	Consider, $A_1$ = $< 105, 104, . . . , 94 >$, $n_1$ = 12, $a_{n_1}$ = 94, $W_1$ = 1194, 
		
		$M$ = 1113 and $W_1 - M$ = 1194-1113 = $81\in\mathbb{N}$. Partition 81 such that 
		
		81 = 1+2+3+4+5+6+7+8+9+10+11+12+13-10+14-(1+2+11) 
		
	\hfill = 3+4+5+6+7+8+9+12+13+14 so that $k_1$ = 10 $<$ 12 = $n_1$, 
		
		$B_1$ = $< 0, 0, 3, 4, 5, 6, 7, 8, 9, 12, 13, 14 >$ and  
		
		$S_1$ = $\{s\in A_1 - B_1\}$ = $\{105-0,104-0, 103-3,102-4,101-5,100-6,$
		
		\hfill $99-7,98-8,97-9,96-12,95-13,94-14\}$ =
		
		~~~ $\{105,104, 100, 98, 96, 94, 92, 90, 88, 84,82,80\}$ with $M_1$ = 1113 = $M$.
		
		\vspace{.1cm}
		\noindent		
		\item [\rm (ii)]	Take 12 (= $n_2$) biggest numbers of $I_{105}\setminus S_1$ as components of $A_2$. Thus, 
		
		$A_2$ = $< 103,102,101, 99, 97, 95, 93, 91, 89, 87,86,85 >$, 
		$n_2$ = 12, 
		
		$a_{n_2}$ = 85, $W_2$ = 1128, $ M$ = 1113 and $W_2 - M$ = 1128-1113 = $15\in\mathbb{N}$.  
			
		\vspace{.1cm}
		15 = 1+2+3+4+5 = 4+5+6 = 2+3+10 are partitions of 15.  
			
		\vspace{.1cm}
		\noindent		
		\item [\rm (a)]	When 15 = 1+2+3+4+5, let $B_2$ = $< 0, 0, 0, 0, 0, 0, 0, 1, 2, 3, 4, 5 >$. 
		
		$\therefore$ $A_2-B_2$  = $< 103, 102, 101, 99, 97, 95, 93,$ 
		
		\hfill $ 91-1, 89-2, 87-3, 86-4, 85-5 >$ 
		
		~~~ ~~   = $< 103, 102,$ $101, 99, 97, 95,$ $93, 90, 87, 84, 82, 80 >$ which is not a subset of $I_{105} \setminus S_1$ and so this partition of 15 is not suitable to consider $S_2$ = $A_2-B_2$. 
			
		\vspace{.1cm}
		\noindent		
		\item [\rm (b)]	When 15 = 4 + 5 + 6, let $B_2$ = $< 0, 0, 0, 0, 0, 0, 0, 0, 0, 4, 5, 6 >$. 
		
		$\therefore$  $A_2-B_2$ = $< 103, 102, 101, 99, 97, 95, 93, 91, 89,$
		
		\hfill $ 87-4, 86-5, 85-6 >$ 
		
		~~ ~~ ~~ ~~ ~~ = $< 103, 102, 101, 99, 97, 95, 93, 91, 89, 83, 81, 79 >$ whose elements are in $\subset I_{105} \setminus S_1$  so that $k_2$ = 3 $<$ 12 = $n_2$, $B_2$ = $< 0, 0, 0, 0, 0, 0, 0,$ $0, 0,$ $4, 5, 6 >$ and  
		 $S_2$ = $\{s\in A_2 - B_2\}$ = $\{ 103, 102, 101, 99, 97, 95, 93, 91, 89,$ $83, 81, 79 \}$ with $M_2$ = 1113 = $M$.
		  	
		 \vspace{.1cm}
		 \noindent		
		 \item [\rm (iii)]	Take 15 (= $n_3$) biggest numbers of $I_{105}\setminus (S_1\cup S_2)$ as components of $A_3$. 
		 
		 Thus, $A_3$ = $< 87, 86, 85, 78,77,76,75,  74,73,72,71,70,69,68,67 >$, $n_3$ = 15,  
		 $a_{n_3}$ = 67, $W_3$ = 1128, $M$ = 1113 and $W_3 - M$ = 1128-1113 = $15\in\mathbb{N}$. 
		 
		 15 = 1+2+3+4+5 is a partition of 15 so that $k_3$ = 5 $<$ 15 = $n_3$. 
		 
		 Let $B_3$ = $< 0, 0, 0, 0, 0, 0, 0, 0, 0, 0, 1, 2, 3, 4, 5 >$.  
		 \\
		 $\therefore$  $A_3-B_3$ = $< 87, 86,85,78,77,76,75,74,73,72,$ 
		 
		 \hfill $71-1,70-2,69-3,68-4,67-5 >$
		 
		~~ ~~ ~~ ~~ ~~  = $< 87, 86, 85,~ 78, 77, 76, 75, 74, 73, 72, ~70, 68, 66, 64, 62 >$~ whose elements are in $\subset I_{105} \setminus (S_1 \cup S_2)$ so that $k_3$ = 5 $<$ 15 = $n_3$, $B_3$ = $< 0, 0, 0, 0, 0, 0, 0, 0, 0, 0, 1, 2,$ $3, 4, 6 >$ and  $S_3 = \{s\in A_3 - B_3\} = \{ 88, 86, 85, 78,77,76,75, 74,73, 72, 70, 68, 66,$ $64, 61 \}$ with $M_3$ = 1113 = $M$.
		 
		 \vspace{.1cm}
		 \noindent		
		 \item [\rm (iv)]	Take 20 (= $n_4$) biggest numbers of $I_{105}\setminus (S_1\cup S_2\cup S_3)$ as components of $A_4$.  
		 
		 $\Rightarrow$	$A_4 = <71,69,67,65,63,62,60,59,58,57,56,55,54,53,$ 
		 
		 \hfill $52,51,50,49,48,4>$, 
		 \\
		 $n_4$ = 20, $a_{n_4}$ = 47, $W_4$ = 1146, $M$ = 1113 and $W_4 - M$ = 1146-1113 =   $33\in\mathbb{N}$.  
		 
		 33 = 1+2+4+5+6+7+8 is a partition of 33 so that $k_4$ = 7 $<$ 20 = $n_4$.   
		 
		 Let $B_4$ = $< 0, 0, 0, 0, 0, 0, 0, 0, 0, 0, 0, 0, 0, 1, 2, 4, 5, 6, 7, 8 >$. 
		 \\
		 $\therefore$  $A_4-B_4$ = $< 71, 69, 67, 65, 63, 62,60,59,58,57,56,55,54,$ 
		 
		 \hfill $53-1,52-2,51-4, 50-5,49-6,48-7,47-8 >$
		 \\
		 whose elements are in  $\subset I_{105} \setminus (S_1 \cup S_2 \cup S_3)$ so that 
		 $k_4$ = 7 $<$ 20 = $n_4$, $B_4$ = $< 0, 0, 0, 0, 0, 0, 0, 0, 0, 0, 0, 0, 0, 1, 2, 4, 5, 6, 7, 8 >$ and $S_4$ = $\{s\in A_4 - B_4\}$ = $< 71, 69, 67, 65, 63, 62,60,59,58,57,56,55,54,52,50,$ 
		 
		 \hfill $47,45,43,41,39 >$ with $M_4$ = 1113 = $M$.
		 
		 \vspace{.1cm}
		 \noindent		
		 \item [\rm (v)]	Now, $S_5$ = $I_{105} \setminus (S_1 \cup S_2 \cup S_3 \cup S_4)$ = $\{53,51,49,48,46, 44,42,40, 38,$ $37,...,1\}$, $n_5$ = 46, $a_{n_5}$ = 1 and $W_5 = M$ = 1113. 
		 
		 Thus, $\{S_1, S_2, S_3, S_4, S_5\}$ is a $5$-partite constant sum partiton of $I_{105}$ with $n_1 = 12 = n_2$, $n_3$ = 15, $n_4$ = 20, $n_5$ = 46 and $M$ = 1113 where 
		 
		 $S_1$ = $\{105,104, 100, 98, 96, 94, 92, 90, 88, 84,82,80\}$,
		 
		 $S_2$ = $\{ 103, 102, 101, 99, 97, 95, 93, 91, 89, 83, 81, 79 \}$,
		 
		$S_3 =  \{ 88, 86, 85, 78,77,76,75, 74,73, 72, 70, 68, 66, 64, 61 \}$,
		
		 $S_4 = \{ 71,69,67,65,63,62,60,59,58,57,56,55,54, 52,$ 
		 
		 \hfill $50,47,45,43,41,39 \}$ and 
		 
		 $S_5$ = $\{53,51,49,48,46, 44,42,40, 38,37,...,1\}$.
		 \\
		  The followings are also 5-partite constant sum partitions of $I_{105}$ with prescribed $S_1$ = $n_1$ = 12 = $n_2$, $n_3$ = 15, $n_4$ = 20, $n_5$ = 46 and $M$ = 1113 using Remark \ref{rem 5}. These are obtained by constant sum partition permutations. Bolded elements are changed elements but still keep their partitions as constant sum partitions. 
		 \begin{enumerate}
		 	\item [\rm (a)] $S_1$ = $\{104, {\bf 103}, {\bf 102}, 98, 96, 94, 92, 90, 88, 84,82,80\}$,
		 	
		 	$S_2$ = $\{ {\bf 105}, 101, {\bf 100}, 99, 97, 95, 93, 91, 89, 83, 81, 79 \}$,
		 	
		 	$S_3 =  \{ 88, 86, 85, 78,77,76,75, 74,73, 72, 70, 68, 66, 64, 61 \}$,
		 	
		 	$S_4 = \{ 71,69,67,65,63,62,60,59,58,57,56,55,54, 52,$ 
		 	
		 	\hfill $50,47,45,43,41,39 \}$,  
		 	
		 	$S_5 = \{53,51,49,48,46, 44,42,40, 38,37,...,1\}$.
		 	
\item [\rm (b)] $S_1 = \{105, 104, 100, 98, 96, 94, 92, 90, 88, 84, 82, 80 \}$,

$S_2$ = $\{ 103, 102, 101, 99, 97, 95, 93, 91, 89, 83, 81, 79 \}$,

$S_3 =  \{ 88, 86, 85, 78,77,76,75, 74,73,72, {\bf 71,} {\bf 67,} 66, 64, 61 \}$,

$S_4 = \{ {\bf 70,} 69, {\bf 68,} 65,63,62,60,59,58,57,56,55,54, 52,$ 

\hfill $50,47,45,43,41,39 \}$,  

$S_5$ = $\{53,51,49,48,46, 44,42,40, 38,37,...,1\}$.

	\item [\rm (c)] $S_1$ = $\{105,104, 100, 98, 96, 94, 92, 90, 88, 84,82,80\}$,
	
	$S_2$ = $\{ 103, 102, 101, 99, 97, 95, 93, 91, 89, 83, 81, 79 \}$,
	
	$S_3 =  \{ 88, 86, 85, 78,77,76,75, 74,73, 72, 70, 68, 66, 64, 61 \}$,
	
	$S_4 = \{ 71, 69, 67,65,63,62,60,59,58,57, 55,54, {\bf 53,} 52,$ 
	
	\hfill $ {\bf 51, 49,} 45,43,41,39 \}$,  

$S_5$ = $\{{\bf 56, 50,} 48, {\bf 47,} 46, 44,42,40, 38,37,...,1\}$.

	\item [\rm (d)] $S_1$ = $\{104, {\bf 103}, {\bf 102}, 98, 96, 94, 92, 90, 88, 84,82,80\}$,
	
	$S_2$ = $\{ {\bf 105}, 101, {\bf 100}, 99, 97, 95, 93, 91, 89, 83, 81, 79 \}$,
	
$S_3 =  \{ 88, 86, 85, 78,77,76,75, 74,73,72, {\bf 71,} {\bf 67,} 66, 64, 61 \}$,

$S_4 = \{ {\bf 70,} 69, {\bf 68,} 65,63,62,60,59,58,57,56,55,54, 52,$

\hfill $50,47,45,43,41,39 \}$,  

$S_5$ = $\{53,51,49,48,46, 44,42,40, 38,37,...,1\}$.

\item [\rm (e)] $S_1$ = $\{104, {\bf 103}, {\bf 102}, 98, 96, 94, 92, 90, 88, 84,82,80\}$,

$S_2$ = $\{ {\bf 105}, 101, {\bf 100}, 99, 97, 95, 93, 91, 89, 83, 81, 79 \}$,

$S_3 =  \{ 88, 86, 85, 78,77,76,75, 74,73,72, {\bf 71,} {\bf 67,} 66, 64, 61 \}$,

$S_4 = \{ {\bf 70,} 69, {\bf 68,} 65,63,62,60,59,58,57,56,$

\hfill $55,54, {\bf 53,} 52, {\bf 51, 49,} 45,43,41,39  \}$,

$S_5$ = $\{{\bf 56, 50,} 48, {\bf 47,} 46, 44,42,40, 38,37,...,1\}$.  
 \end{enumerate} 
One can find more such 5-partite constant sum partitions of $I_{105}$ with prescribed $S_1$ = $n_1$ = 12 = $n_2$, $n_3$ = 15, $n_4$ = 20, $n_5$ = 46 and $M$ = 1113. \hfill $\Box$		 } 
\end{exm}

\section{Applications of Constant Sum Partitions of $I_n$ in Graph theory}

Partitioning of $I_n$ has applications in graph theory. A graph $G$ with $n(n+1)/2$ edges is called {\it a path-perfect graph} if the edge set of $G$ can be partitioned as $E_1 \cup E_2 \cup . . . \cup E_n$ so that $E_i$ induces a path of length $i$, for $1 \leq i \leq n$.  In \cite{ss}, existence of $(s,t)$-partition of $I_n$ is established whenever $n(n+1)$ = $2st$ and also an $r$-regular path-perfect graph $G$ with $p$ vertices and $n(n+1)/2$ edges using the existence of an $(r/2, p)$-partition of $I_n$ when $r$ is even. 

Alavi \cite{ab} introduced the concept of {\it Ascending Subgraph Decomposition (ASD)} of a graph $G$ with size $\binom{n+1}{2}$ as the decomposition of $G$ into $n$ subgraphs $G_1,G_2,\ldots,G_n$ without isolated vertices such that each $G_i$ is isomorphic to a proper subgraph of $G_{i+1}$  and  $|E(G_i)|  = i$ for $1 \leq i\leq n$. Nagarajan \cite{nn} generalized ASD to $(a, d)$-ASD of graph $G$ with size $\binom{(2a+(n-1)d)n}{2}$ as the decomposition of $G$ into $n$ subgraphs $G_1,G_2,\ldots,G_n$ without isolated vertices such that each $G_i$ is isomorphic to a proper subgraph of $ G_{i+1}$ and $|E(G_i)|  = a+(i-1)d$ for $1 \leq i \leq n$. Gnana Dhas \cite{gp} defined {\it $(a, d)$-Continuous Monotonic Decomposition} or $(a, d)$-CMD of a graph $G$ of size $\binom{(2a+(n-1)d)n}{2}$ as the decomposition of $G$ into $n$ subgraphs $G_1,G_2,\ldots,G_n$ such that each $G_i$ is connected and  $|E(G_i)| = a+(i-1)d$ for $i = 1,2,\ldots,n$. Clearly, CMD of a graph $G$ and its $(1,1)$-CMD are the same.

Vilfred \cite{vs} defined {\em $(a,d)$-Continuous Monotonic Subgraph Decomposition} or $(a,d)$-CMSD of graphs as graphs which have both $(a,d)$-ASD and $(a,d)$-CMD and developed its theory, $a,d\in\mathbb{N}$. Clearly, $(a,d)$-CMSD of graph $G$ with size $\binom{(2a+(n-1)d)n}{2}$ is the decomposition of $G$ into $n$ subgraphs $G_1,G_2,\ldots,G_n$ without isolated vertices such that each $G_i$ is {\it connected and isomorphic to a proper subgraph} of $G_{i+1}$ and  $|E(G_i)| = a+(i-1)d$ for $1\leq  i\leq n$.

Vilfred \cite{v22} defined {\em $k$-distance magic labeling ($k$-DML)} of a graph $G$ of order $n \geq 3$ as a labeling $f:$ $V(G)$ $\rightarrow$ $\{1,2,...,n\}$ such that $\sum_{w\in \partial N_k(u)}{ f(w)}$ is a constant and independent of $u\in V(G)$ where $\partial N_k(u)$ = $\{v\in V(G):$ $d(u, v) = k\}$, $k\in\mathbb{N}$. Graph $G$ is called a \emph{$k$-distance magic ($k$-DM)} if it has a $k$-DML. $k$-DM labeling is a generalisation of DM labeling or Sigma labeling  \cite{v96,v22} of graphs, $k\in\mathbb{N}$. Partitioning of integers and these decomposition of graphs are very much related. For further reading on applications of partitioning of $I_n$ in graph theory, one can refer \cite{ab,gp,mz,ss,v96,v22,vs}.

\end{document}